\documentclass[a4paper,10pt]{article}

\usepackage[T2A]{fontenc} 
\usepackage[utf8]{inputenc}
\usepackage{amsmath}
\usepackage{amssymb}
\usepackage[english,russian]{babel}
\usepackage{graphicx}

\newcommand{\eps}{\varepsilon}

\begin{document}
\large

{\it ДОКЛАДЫ АКАДЕМИИ НАУК, 201X, том X, \No Y.}

\begin{center}{\bf = МАТЕМАТИКА =}\end{center}

\large
\begin{center}
\it \quad УДК 517.925.42 \hfill \quad  \ \\
 \bf ВЫРОЖДЕННЫЕ РЕЗОНАНСЫ И ИХ УСТОЙЧИВОСТЬ В ДВУМЕРНЫХ СИСТЕМАХ С МАЛОЙ ОТРИЦАТЕЛЬНОЙ ДИВЕРГЕНЦИЕЙ\\
\end{center}
\centerline{ \bf  О. Ю. Макаренков, И. С. Мартынова}
\vskip0.3cm \centerline{  \normalsize Представлено академиком С. Н. Васильевым} \vskip0.5cm

\renewcommand{\figurename}{Фиг.}
\large

\setcounter{equation}{0} \thispagestyle{empty}

В работе показывается, что связанные с именами Боголюбова \cite{bog}, Малкина \cite{mal} и Мельникова \cite{mel} классические условия бифуркации $T$-периодических решений в аналитических дифференциальных уравнениях вида
\begin{equation}\label{ps}
  \dot x =f(x)+\eps g(t,x,\eps),\quad x\in\mathbb{R}^2
\end{equation}
в окрестности $T$-периодического цикла $x_0$ соответствующей порождающей системы
\begin{equation}\label{np}
  \dot x=f(x),\quad x\in\mathbb{R}^2
\end{equation}
могут быть существенно ослаблены, коль скоро возмущенная система (\ref{ps}) обладает  следующим свойством отрицательности дивергенции
\begin{equation}\label{div}
  \sum_{i=1,2} (f_i)'_{x_i}+\eps (g_i)'_{x_i}(t,x,\eps)<0,\quad t\in\mathbb{R},\ x\in V,\ \eps\in(0,\eps_0).
\end{equation}
Здесь $V\in\mathbb{R}^2$ -- какая-нибудь окрестность кривой $t\mapsto x_0(t)$  и $\eps_0>0$ -- любая подходящая константа.

Результаты классической теории возмущений \cite{bog,mal,mel} связаны с разложением оператора Пуанкаре $\mathcal{P}_\eps$ за период $T>0$ возмущенной системы (\ref{ps}) по степеням малого параметра $\eps$:  
$$
   \mathcal{P}_\eps(\phi(x,\eps))=\phi(x,\eps)+\eps^\lambda \overline{f}(x)+\eps^{2\lambda} \overline{\overline{f}}(x)+...,
$$
где $\lambda\in(0,1]$, $\phi(\cdot,\eps)$ -- подходящая взаимооднозначная замена переменных и $\overline{f},$ $\overline{\overline{f}}$ и т. д.  -- так называемые бифуркационные функции 1-го, 2-го и т. д. порядков.
Для бифуркации $T$-периодического решения системы (\ref{ps}) из решения $x_0$, таким образом, необходимо, чтобы $x_0(0)$ был нулем бифуркационной функции первого порядка. Невырожденность данного нуля, то есть обратимость производной бифуркационной функции первого порядка в $x_0(0)$, является достаточным условием бифуркации. Невырожденность нуля бифуркационной функции первого порядка имеет место в следующих основных ситуациях
\begin{itemize}
\item[1.] Боголюбов \cite{bog}: все решения системы (2) -- $T$-периодические и $x_0(0)$ является простым нулем функции усреднения $$\overline{f}(v)=\int_0^T \Omega'_x(0,\tau,\Omega(\tau,0,v))g(\tau,\Omega(\tau,0,v),0)d\tau,$$ где $\Omega(\cdot,t_0,v)$ -- решение системы (\ref{np}) с начальным условием $x(t_0)=v.$ 

\item[2.] Мельников \cite{mel}: решение системы (2) с любым начальным условием $v$ -- периодическое с периодом $T(v)$, $s=0$ является простым нулем период-функции $s\mapsto T(x_0(0)+\dot x_0(0)^\perp s),$ где $v^\perp={\rm col}(-v_2,v_1),$ и $\theta=0$ является простым нулем субгармонической бифуркационной функции Мельникова $\theta\mapsto\pi_E(x_0(\theta)) \overline{f}(x_0(\theta)),$ где $\pi_E(x_0(\theta))$ -- определенная 2x2-матрица.

\item[3.] Малкин \cite[Гл. VI, \S 2]{mal}: матрица $\Omega'_x(T,0,x_0(0))$ имеет отличное от $\pm 1$ собственное значение и так называемая бифуркационная функция Малкина-Луда $\theta\mapsto\pi_A(x_0(\theta)) \overline{f}(x_0(\theta)),$ где $\pi_A(x_0(\theta))$ -- 2x2-матрица, имеет $\theta=0$ простым нулем.
\end{itemize}

Если ни одно из указанных 3-х типов условий не выполнено, приходится задействовать бифуркационную функцию $\overline{\overline{f}}$ 2-го и высших порядков. Полученные периодические решения называются в этом случае вырожденными резонансами. Для доказательства существования вырожденных резонансов оказывается полезной имеющаяся в системе дополнительная структура, см. Тхай \cite{Tkhai}. В общем же случае условия бифуркации выписываются в терминах решений вспомогательных полиномиальных уравнений, см. Копнин \cite{kop}, Yagasaki \cite{yag}. 

В работах многих авторов (см. ссылки в \cite{maw,kmn,mak}) показывалось, что условия на невырожденность нулей бифуркационных функций может быть заменено невырожденностью их топологического индекса. Однако устойчивость соответствующих периодических решений возмущенной системы явно доказана (Макаренков-Ортега \cite{mak_ort}) только для случая Малкина. В этой статье, пользуясь одним лишь топологическим индексом и ограничиваясь случаем возмущенных систем (\ref{ps}) с аналитическими правыми частями и отрицательной дивергенцией (быть может не имеющей места при $\eps=0$), предлагаются результаты об  устойчивости также в случаях Боголюбова и Мельникова.

1. Рассмотрим сначала ситуацию Боголюбова. Напомним, что система (\ref{ps}) с $T$-периодической по $t$ правой частью является аналитической в некоторой окрестности, если правая часть разлагается в этой окрестности в ряд по степеням фазовой переменной и сходимость рассматриваемого ряда равномерна по $t$ и $\eps.$ Ниже мы используем понятие индекса Пуанкаре ${\rm ind}(v_0,\overline{f})$ изолированного нуля $v_0$ векторного поля $\overline{f},$ который определяется как топологическая степень ${\rm d}(\overline{f},U)$ поля $\overline{f}$ относительно границы малой отрытой окрестности $U$ точки $v_0$, см. \cite{kra_zab}. Топологическую степень в  рассматриваемой двумерной ситуации можно определить как число полных оборотов, совершаемых вектором $\overline{f}(v(\theta))$, в то время как $v(\theta)$ обходит $v_0$, при изменении $\theta$ от $0$ до $2\pi$, против часовой стрелки вдоль границы $\partial U$ окрестности $U$. При этом, каждый из таких полных оборотов засчитывается с положительным или отрицательным знаком в зависимости от того достигает вектор $\overline{f}(v(\theta))$ вектора $\overline{f}(v(0))$ против часовой стрелки или по ее направлению, см. \cite{kra_zab}.

{\bf Теорема 1.} {\it  Предположим, что рассматриваемая возмущенная система (\ref{ps}) имеет $T$-периодические по $t$ правые части и является аналитической в малой окрестности $V$ $T$-периодического решения $x_0$, а все решения порождающей системы (\ref{np}) из данной окрестности являются $T$-периодическими. Если выполнено условие отрицательности дивергенции (\ref{div}) и  
\begin{equation}\label{index}
    {\rm ind}(x_0(0),-\overline{f})>0,
\end{equation}
то существует $\eps_0>0$ такое, что при $\eps\in(0,\eps_0]$ система (\ref{ps}) имеет по крайней мере одно асимптотически устойчивое $T$-периодическое решение $x_\eps$ такое, что $x_\eps\to x_0$ при $\eps\to 0.$
}

{\bf Пример 1.} Рассмотрим следующий вариант уравнения Дуффинга:   $\ddot u+\eps c\dot u+\eps^3au+\eps^2 b u^3=\eps^2\gamma\cos\omega t$, в котором $a,b,c,\gamma,\eps>0.$ Функция $u$ является решением данного уравнения тогда и только тогда, когда $x=(u,(1/\eps)\dot u)$ -- решение системы
\begin{equation}\label{duf}
\begin{array}{lll}
  \dot x_1 & = & \eps x_2,\\
  \dot x_2 & = & -\eps c x_2-\eps^2 a x_1-\eps b (x_1)^3+\eps\gamma\cos\omega t. 
\end{array}
\end{equation}
Соответствующая функция усреднения $\overline{f}$ выписывается как 
$
  \overline{f}(v)={\rm col}(v_2,-b v_1^3-cv_2),
$
и имеет единственный нуль $v_0=0.$ Поэтому, бифуркация $2\pi/\omega$-периодического решения системы (\ref{duf}) возможна из одного лишь решения $x_0(t)\equiv 0.$ Однако, ${\rm det}\|\overline{f}'(0)\|=0$ и усреднение по первому приближению не позволяет доказать существование бифуркации. Вместо вычисления высших приближений применим теорему 1, поскольку вычисление дивергенции в формуле (\ref{div}) приводит для системы (\ref{duf}) к выражению $-\eps c,$ то есть является отрицательной. Замечая, что векторное поле $\overline{f}$ линейно гомотопно тождественному на мылых окружностях с центрами в $v_0=0,$ приходим к заключению ${\rm ind}(0,-\overline{f})=1$. Значит, при всех достаточно малых $\eps>0$ система (\ref{duf}) действительно имеет по крайней мере одно асимптотически устойчивое $2\pi/\omega$-периодическое решение, сходящееся к 0 при $\eps\to 0.$

2. Перейдем к ситуации Мельникова. Ниже, под сходимостью векторов $v_\eps\in\mathbb{R}^2$ к множеству $x_0(\mathbb{R})$ при $\eps\to 0$ понимается сходимость в метрике Хаусдорфа.

{\bf Теорема 2.} {\it  Предположим, что рассматриваемая возмущенная система (\ref{ps}) имеет $T$-периодические по $t$ правые части и является аналитической в малой окрестности $V\subset\mathbb{R}^2$ $T$-периодического цикла $x_0$. Пусть далее ни одно из решений порождающей системы (\ref{np}) из множества $V\backslash x_0(\mathbb{R})$ 
не является $T$-периодическим. Если выполнено условие отрицательности дивергенции (\ref{div}) и  
\begin{equation}\label{index2}
    {\rm d}(-\overline{f},U)\not=1,
\end{equation}
где $U$ -- внутренность цикла $x_0$,
то существует $\eps_0>0$ такое, что при $\eps\in(0,\eps_0]$ система (\ref{ps}) имеет по крайней мере одно асимптотически устойчивое $T$-периодическое решение $x_\eps$ такое, что $x_\eps(0)\to x_0(\mathbb{R})$ при $\eps\to 0.$ Начальное условие $x_\eps(0)$ приближается к $x_0(\mathbb{R})$ изнутри или снаружи в зависимости от того ${\rm d}(-\overline{f},U)>1$ или ${\rm d}(-\overline{f},U)<1.$}

Проверка условия (\ref{index2}) может проводиться на основании анализа классических бифуркационных функций Малкина и Мельникова. Приведем один такой критерий \cite{mak}. Обозначим через $z_A$ решение уравнения $\dot z=-(f'(x_0(t)))^* z$ с начальным условием $z_A(0)=\left(1/\|\dot{x}_0(0)\|^2\right)\dot{x}_0(0)$, а через $z_E$  -- решение данной линейной системы с начальным условием $\widetilde z_E(0)={\rm col}(-\dot{x}_{0,2}(0),\dot{{x}}_{0,1}(0)).$ Положим
$$
   M_{j}(\theta)=\int_\theta^{\theta+T}\left<z_{j}(\tau),g(\tau-\theta,x_0(\tau),0)\right>d\tau,\quad j=A,E.
$$
Если $M_E$ имеет на $[0,T)$ ровно два нуля $\theta_1,\theta_2$ и является в них строго монотонной, а для функции $M_A$ выполнено условие $M_A(\theta_1)M_A(\theta_2)<0,$ то ${\rm d}(-\overline{f},U)$ принимает либо значение 0, либо значение 2.

{\bf Пример 2.} Рассмотрим систему
\begin{equation}\label{yag}
  \begin{array}{lll}
    \dot x_1&=&x_2\left(\dfrac{1}{4}(x_1^2+x_2^2-2)^p+1\right)\\
    \dot
    x_2&=&-x_1\left(\dfrac{1}{4}(x_1^2+x_2^2-2)^p+1\right)-\eps^2 x_2+\varepsilon\sin(t).
  \end{array}
\end{equation}
Порождающая система допускает семейство циклов
$
  {x}_{0,\alpha}(t)={\rm col}\left(
  \alpha\sin\left(({2\pi}/{T(\alpha)})t\right),
  \alpha\cos\left(({2\pi}/{T(\alpha)})t\right)
  \right)
$
 с периодами
 $
   T(\alpha)={2\pi}/{\left((1/4)(\alpha^2-2)^p+1\right)}.
 $
При $\alpha=\sqrt{2}$ для период-функции имеем 
 $
   T\left(\alpha\right)=2\pi,$ $T'\left(\alpha\right)=...=T^{(p-1)}\left(\alpha\right)=0,$
 $T^{(p)}\left(\alpha\right)\not=0$ и усреднение по первому приближению (или стандартный метод Мельникова) не применимо. Вычисление высших приближений для данного примера проведено в \cite{yag}. Для дивергенции системы (\ref{yag}) получаем значение $-\eps^2$, то есть условие (\ref{div}) выполнено. Поэтому, вместо привлечения высших приближений мы используем теорему 2. При любом $p\in\mathbb{N}$ для функций $M_A$ и $M_E$ получаем выражения 
$ M_A(\theta)=-\frac{1}{\sqrt{2}}\pi\cos(\theta),$ $
  M_E(\theta)=-\sqrt{2}\pi\sin(\theta).$ 
Условия указанного выше критерия выполнены (два нуля на интервале $[0,2\pi)$ даются формулами $\theta_1=0$ и
$\theta_2=\pi$), поэтому условие (\ref{index2}) выполнено. Применяя теорему 2, заключаем, что при всех достаточно малых $\eps>0$ система (\ref{yag}) имеет асимптотически устойчивое $2\pi$-периодическое решение, сходящееся при $\eps\to 0$ к кругу радиуса $1/\sqrt{2}.$

3. Обсудим основные идеи доказательства предложенных теорем. Для получения результатов теорем 1 и 2 прежде всего доказывается положительность топологической степени оператора Пуанкаре возмущенной системы относительно окрестности точки $x_0(0)$ в случае теоремы 1  и относительно внутренней и внешней окрестностей множества $x_0(\mathbb{R})$ в случае теоремы 2. В первом случае используется результат Мавена \cite{maw}, а во втором теорема Каменского-Макаренкова-Нистри \cite{kmn} в сочетании с результатом Capietto-Mawhin-Zanolin \cite{maw1}. Далее используется свойство (\ref{div}) и теорема из Nakajima-Seifert \cite{nak}, которая позволяет заключить существование у возмущенной системы (\ref{ps}) изолированного $T$-периодического решения с индексом Пуанкаре +1. Наконец,
асимптотическая устойчивость такого решения следует из теоремы Колесова-Ортега \cite{kol,ort}. 

Исследования авторов поддержаны грантом Президента РФ для молодых кандидатов наук (МК-1530.2010.1) и грантом РФФИ 09-01-00468. Исследования первого автора также поддержаны грантом РФФИ 10-01-93112.


\

\noindent {\bf О.~Ю.~Макаренков (автор для
переписки):} Институт проблем управления им. В.А. Трапезникова, 117997, ул. Профсоюзная 65, Москва, Россия и Имперский Колледж Лондона, SW7 2AZ, Великобритания ({omakarenkov@mail.ru}).\\
Служебный телефон: (495) 334-86-60

\vskip0.2cm \noindent {\bf И.~С.~Мартынова:} Воронежская государственная технологическая академия, 394036, пр. Революции 19, Воронеж, Россия ({i\_martynova@inbox.ru}).\\
Служебный телефон: (4732) 552-550

\end{document}